\documentclass[a4paper,12pt,reqno]{amsart}
\usepackage[french,english]{babel}
\usepackage{amsmath, amsfonts, amssymb, amsthm,mathptmx}
\usepackage{color}
\usepackage{graphicx}
\usepackage{multirow}

\allowdisplaybreaks

\usepackage[letterpaper, hmargin=1in, top=1in, bottom=1.2in, footskip=0.7in]{geometry}

\usepackage{rotating,pdflscape,color,amssymb,subfigure,psfrag,amsmath,eufrak,bbm,epsfig}
\definecolor{vdarkred}{rgb}{0.6,0,0.2}
\definecolor{vdarkblue}{rgb}{0,0.2,0.6}
\usepackage[pdftex, colorlinks, linkcolor=vdarkblue,citecolor=vdarkred]{hyperref}

\newtheorem{lemme}{Lemma}[section]

\newtheorem{theo}[lemme]{Theorem}

\newtheorem{cor}[lemme]{Corollary}

\newtheorem{propo}[lemme]{Proposition}

\newenvironment{dem*}{\underline{\textbf{Proof }}}{\hfill$\square$ $$$$}

\theoremstyle{definition}
\newtheorem{assum}{Hypothesis} 
\newtheorem{defi}[lemme]{Definition}
\newtheorem{rem}[lemme]{Remark}

\DeclareMathOperator{\im}{\operatorname{Im}}
\DeclareMathOperator{\re}{\operatorname{Re}}

\DeclareMathOperator{\cloi}{\overset{(d)}{\longrightarrow}}

\DeclareMathOperator{\E}{\mathbb{E}}

\DeclareMathOperator{\pro}{\mathbb{P}}

\DeclareMathOperator{\var}{\mathbb{V}ar}
\DeclareMathOperator{\op}{op}

\DeclareMathOperator{\diag}{\text{diag}}

\DeclareMathOperator{\C}{\mathbb{C}}
\DeclareMathOperator{\R}{\mathbb{R}}

\newcommand{\ope}{\operatorname}
\newcommand{\tto}{\longrightarrow}
\newcommand{\ol}{\overline}
\newcommand{\wt}{\widetilde}

\newcommand{\ep}{\varepsilon}

\newcommand{\wg}{\operatorname{Wg}}

\newcommand{\tr}{\operatorname{Tr}}
\newcommand{\dI}{\mathrm{I}}
\newcommand{\dV}{\mathrm{V}}
\newcommand{\bA}{\mathbf{A}}

\newcommand{\bV}{\mathbf{V}}
\newcommand{\bU}{\mathbf{U}}

\newcommand{\bP}{\mathbf{P}}
\newcommand{\bB}{\mathbf{B}}

\newcommand{\bI}{\mathbf{I}}

\newcommand{\bQ}{\mathbf{Q}}
\newcommand{\bR}{\mathbf{R}}

\newcommand{\be}{\mathbf{e}}

\newcommand{\bM}{\mathbf{M}}
\newcommand{\bN}{\mathbf{N}}
\newcommand{\bX}{\mathbf{X}}

\newcommand{\bW}{\mathbf{W}}

\newcommand{\bH}{\mathbf{H}}
\newcommand{\bE}{\mathbf{E}}
\newcommand{\bJ}{\mathbf{J}}
\newcommand{\bY}{\mathbf{Y}}

\newcommand{\tred}{\textcolor{red}}

\newcommand{\bpr}{\begin{pr}}
\newcommand{\epr}{\end{pr}}
\newcommand{\trm}{\textrm}
\newcommand{\f}{\frac}

\newcommand{\ff}{\frac{1}}
\newcommand{\one}{\mathbbm{1}}

\newcommand{\ds}{\displaystyle}
\newcommand{\lf}{\left}
\newcommand{\ri}{\right}
\newcommand{\st}{such that }

\newcommand{\Tr}{\operatorname{Tr}}

\newcommand{\bes}{\begin{equation*}}
\newcommand{\ees}{\end{equation*}}
\newcommand{\beqy}{\begin{eqnarray}}
\newcommand{\eeqy}{\end{eqnarray}}
\newcommand{\beq}{\begin{eqnarray*}}
\newcommand{\eeq}{\end{eqnarray*}}
\newcommand{\bbe}{\begin{equation}}
\newcommand{\ee}{\end{equation}}
\newcommand{\bbm}{\begin{bmatrix}}
\newcommand{\ebm}{\end{bmatrix}}
\newcommand{\bpm}{\begin{pmatrix}}
\newcommand{\epm}{\end{pmatrix}}
\newcommand{\bdet}{\begin{vmatrix}}
\newcommand{\edet}{\end{vmatrix}}
\newcommand{\la}{\label}
\newcommand{\eqre}{\eqref}
\newcommand{\ti}{\times}

\newcommand{\egd}{\ := \ }

\newcommand{\OO}[1]{O \left( #1 \right)}
\newcommand{\oo}[1]{o \left( #1 \right)}

\newcommand{\ie}{i.e. } 
\newcommand{\id}{\ope{id}} 
\newcommand{\ovl}{\overline}

\newcommand{\alp}{\alpha}
\newcommand{\al}{\alp}

\newcommand{\bet}{\beta}
\newcommand{\ga}{\gamma}
\newcommand{\tta}{\theta}
\newcommand{\lam}{\lambda}
\newcommand{\si}{\sigma}

\newcommand{\ld}{\ldots}
\newcommand{\norm}[1]{\left\| #1 \right\|_{\op}}

\newcommand{\lexp}[2]{{\vphantom{#2}}^{#1}#2}
\newcommand{\lbinom}[3]{{\vphantom{#3}}^{\;\;#1}_{#2}#3}

\newcommand{\Ec}[1]{\E \left[ #1 \right]}
\newcommand{\inve}[1]{ \left( #1 \right)^{-1}}

\newcommand{\bgt}{\begin{itemize}}
\newcommand{\ent}{\end{itemize}}
\newcommand{\ite}{\item}

\newenvironment{pr}{\noindent {\bf Proof. }}{\hfill $\square$\\}

\newcommand{\ii}{\operatorname{i}}

\newcommand{\perm}{\mathfrak{S}}
\newcommand{\fix}{\operatorname{Fix}}
\newcommand{\Cat}{\operatorname{Cat}}

\title{Fluctuation of matrix entries and application to outliers of elliptic matrices}
\author{Florent Benaych-Georges, Guillaume C\'ebron and Jean Rochet}

\thanks{FBG and JR: MAP5,
Universit\'e Paris Descartes,
45, rue des Saints-P\`eres
75270 Paris Cedex 06, France. florent.benaych-georges@parisdescartes.fr, jean.rochet@parisdescartes.fr.\\ \indent GC: IMT,
Universit\'e Paul Sabatier,
118 Route de Narbonne
31062 Toulouse Cedex 04, France. guillaume.cebron@math.univ-toulouse.fr. (GC was partially supported by the ERC advanced grant ``non-commutative distributions in free probability", held by R. Speicher).}



\begin{document}
\maketitle
\begin{abstract}
For any family of $N\times N$ random matrices $(\mathbf{A}_k)_{k\in K}$ which is invariant, in law,  under unitary conjugation, we give general sufficient conditions for central limit theorems for random variables of the type $\operatorname{Tr}(\mathbf{A}_k \mathbf{M})$, where  the matrix  $\mathbf{M}$ is deterministic  (such random variables include for example the  normalized matrix entries of the $\mathbf{A}_k$'s). A consequence  is the asymptotic independence of the projection of the matrices $\mathbf{A}_k$ onto the subspace of null trace matrices from  their projections onto the orthogonal of this subspace. 
These results are used to study the asymptotic behavior of the outliers of a spiked elliptic random matrix. More precisely, we show that the fluctuations of these outliers around their limits can have various rates of convergence, depending on the Jordan Canonical Form of the additive perturbation. Also, some correlations can arise between outliers at a macroscopic distance from each other. These phenomena have already been observed in \cite{FloJean} with random matrices from the Single Ring Theorem. 
\end{abstract}

\section{Introduction}

This paper is first concerned with the fluctuations of linear functions   of entries of  unitarily invariant random matrices   when the dimension tends to infinity. Then, it deals with  the  application of such limit theorems to the fluctuations of the outliers of spiked elliptic matrices.  

The first problem 
is to find out conditions under which, for given collections $(\bA_k)_{k\in K}$ of random matrices and $(\bM_\ell)_{\ell\in L}$ of non-random matrices, the finite marginals of  \bbe\la{17101516h0Intro}
  \left(\tr \big( \bA_{k} \bM_{\ell}  \big)  -\E \tr \big( \bA_{k} \bM_{\ell}  \big)\ri)_{k\in K, \ell\in L}\end{equation} converge as the dimension $N$ tends to infinity. 
  We shall always suppose that   the $ \bA_k$'s and the $\bM_\ell$'s   have Euclidean norms of order $\sqrt{N}$, \ie that the random variables   $$\ff{N}\Tr \bA_k \bA_k^*\qquad\trm{ and }\qquad \ff{N}\Tr \bM_\ell \bM_\ell^*$$ are bounded in probability.
  The case \bbe\la{17101516h0IntroIntro}\bM_\ell=\sqrt{N}\times\text{(an elementary matrix)}\ee is a classical example.    In this framework, 
  the main hypothesis we need for 
the random vector of \eqre{17101516h0Intro} to be asymptotically Gaussian is the global invariance, in law, of $(\bA_k)_{k\in K}$ under unitary conjugation, \ie that for    any non random unitary matrix $\bU$,  $$(\bA_k)_{k\in K} \stackrel{\text{law}}{=}\big(\bU\bA_k \bU^*\big)_{k \in K}.$$
It then appears that the question decomposes into two independent problems: one associated to the projections of the $\bA_k$'s onto the space of null trace matrices (see Theorem \ref{Th_Mainbis}) and one associated to the convergence of 
the centered traces of the $\bA_k$'s;  and that both give rise to independent asymptotic fluctuations (see Theorem \ref{Th_Main} and Corollary \ref{Th_Mainter115}). These results extend an already proved partial result in this direction,  Theorem 6.4 of \cite{BEN2} (see also Theorem 1.2 of \cite{Rains98} in the particular case of real symmetric matrices $\bA_k$). The main advantages of  Theorems \ref{Th_Mainbis} and \ref{Th_Main}  over the results of \cite{Rains98} and \cite{BEN2} is, firstly,  that they do not require   the matrices $\bM_\ell$ to  have  finitely many non zero entries (or to be well approximated by such matrices) and, secondly,  that they give the asymptotic independence mentioned above. Besides, the technical hypotheses needed here are weaker than in the existing literature. Our proofs are based on the so-called \emph{Weingarten calculus}, an integration method for the Haar measure on the unitary group developed by Collins and \'Sniady in \cite{collinsIMRN,COL}.

All these results belong to a long list of theorems begun in 1906 with the theorem  by Borel \cite{b1906} stating that any coordinate of a uniformly distributed random vector of the sphere of $\R^N$ with radius $\sqrt{N}$ is asymptotically standard Gaussian as $N\to\infty$, and continued e.g. with the papers  \cite{Rains98,diaconis2003,jiang06,meckes08, meckessourav08,collins-stolz08,BENCLT,ToddGuillaume,FloJean2} on central limit theorems on large matrix   spaces. Some of the results from these papers can be deduced from this paper (see e.g. Remark \ref{list_applications_240217}). 

Second order freeness,  a  theory that has been developed these last ten years, deals  with Gaussian  fluctuations (called {\it second order limits}) of traces of large random matrices. The most emblematic articles in this theory   are \cite{mingo-nica04,mingo-speicher06, Mingo2007, mingo-piotr-collins-speicher07}. As explained in Remark \ref{no_second_order_freeness_supposed}, our results cannot be deduced   from this theory, because the ``test matrices" we consider (\ie the matrices $\bM_\ell$) are not supposed to have second order limits. Precisely,   in   classical  applications of our results (e.g.  the case  of \eqref{17101516h0IntroIntro}), the matrices $\bM_\ell$ do not have any second order limit.  However, we shall see in Section \ref{21021720h} that    our results    extend the consequences  of the existence of a second order limit   for  unitarily invariant matrix ensembles.

\begin{figure}[ht]
\centering
{
\includegraphics[scale=0.51]{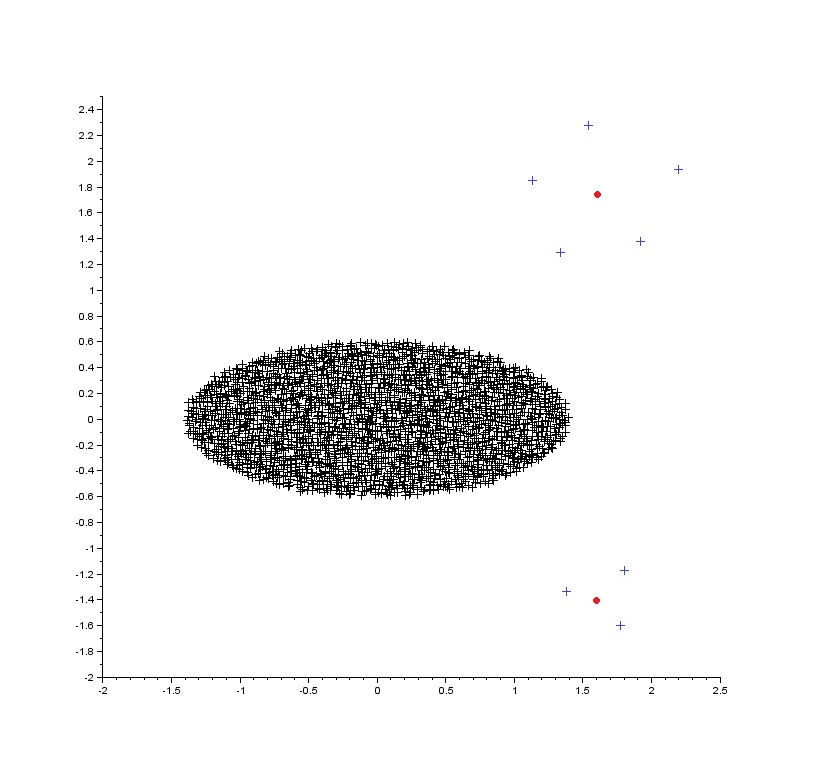}
}
\caption{Spectrum of a Gaussian elliptic matrix of size $N=2500$ with perturbation matrix $\bP = \diag\lf( \bR_5(1.5+2.625\ii),\bR_3(1.5-1.5\ii),0,\ld,0\ri)$. We   see the blue crosses ``\textcolor{blue}{$+$}'' (outliers) forming respectively a regular pentagon and an equilateral triangle around the red dots ``\tred{$\bullet$}'' (their limit).  We also see a significant difference between the two rates of convergence, $N^{-1/10}$ and $N^{-1/6}$.}\label{ellipsecolor}
\end{figure}

The general results about asymptotic fluctuations of matrix entries that we prove here are then  applied  to the fluctuations of the outliers of Gaussian elliptic matrices. 
One can prove (see e.g. \cite{DjaJTP}) that the global behavior of the spectrum  of  a large random matrix is not altered, from the macroscopic point of view, by a low rank additive perturbation. However, some of the eigenvalues, called \emph{outliers}, can deviate away from the bulk, depending on the strength of the perturbation. Firstly brought to light for empirical covariance matrices by Johnstone in \cite{john:01}, this phenomenon,   known as the \emph{BBP transition}, was proved by Baik, Ben Arous and P\'ech\'e in \cite{BBP05}, and then extended  to  several   Hermitian models in  \cite{SP06,FP07,CAP2,CAP,PRS11, BEN1,BEN4,BEN2,BEN3,CAP3,KYIsoSCL,KYIOLW}. Non-Hermitian models have been also studied: i.i.d. matrices in \cite{TAO,BORCAP1,RAJ}, real elliptic matrices in \cite{or13}, matrices from the Single Ring Theorem in \cite{FloJean} and also nearly Hermitian matrices \cite{Rochet,or15}. 
 As an application of our main result, we investigate the fluctuations of the outliers and due to the non-Hermitian structure, we prove, as in \cite{FloJean,BORCAP1, RAJ,Rochet}, that the distribution of the fluctuations highly
depends on the shape of the Jordan Canonical Form of the perturbation. In particular, the convergence
rates depend on the sizes of the Jordan blocks. Also, the outliers tend to locate around their limits at the
vertices of regular polygons (see Figure \ref{ellipsecolor}). At last, as in \cite{FloJean}, we prove the quite surprising fact that outliers at
 macroscopic distance from  each other can have correlations fluctuations  (see Remark \ref{rem1519kl2901}), . 

The paper is organized as follows. In Section~\ref{mainresults}, we state our main results (Theorems~\ref{Th_Mainbis}, \ref{Th_Main}, \ref{theo08521412020154} and   \ref{thbelow0911141210002015}) and their corollaries. These theorems are then proved in  the following  sections and an appendix is devoted to a technical result needed here.

\noindent{\bf Notation:} For $u,v$ sequences, $u=o(v)$ means that $u/v\to 0$ and $u=\OO{v}$ means that $u/v$ is bounded. Also, the dimension $N$ of the matrices is most times an implicit parameter.
\section{Main results}\label{mainresults}

\subsection{General results} Let $\bA  = \big(\bA_k\big)_{k \in K}$ be a collection of   $N\ti N$ random  matrices and let $  \big(\bM_\ell\big)_{\ell \in L}$ be a collection $N\ti N$ non random   matrices, both implicitly depending on $N$. 
\begin{assum}\la{Assum201611}\noindent\bgt
\item[(a)]  $\bA $ is  invariant in distribution under  unitary conjugation: for    any non random unitary matrix $\bU$,  $$\big(\bA_k\big)_{k \in K} \stackrel{\text{law}}{=}\big(\bU\bA_k \bU^*\big)_{k \in K};$$
\item[(b)]  for each $k\in K$, and each $p,q\geq 1$,  $\ff{N} \tr (\bA_k\bA_k^*)^{p}$ is bounded in $L^q$ independently of $N$;
\item[(c)] for each $k,k'\in K$, we have the following convergences, in $L^2$,  to deterministic limits: 
$$\lim_{N \to \infty}\ff{N} \tr \bA_k\bA_{k'}-\ff{N} \tr \bA_k\cdot \ff{N} \tr\bA_{k'}=\tau(k,k')$$and $$\lim_{N \to \infty}\ff{N} \tr \bA_k\bA_{k'}^*-\ff{N} \tr \bA_k\cdot \ff{N} \tr\bA_{k'}^*=\tau(k,\ovl{k'});$$
\item[(d)] for each $\ell,\ell'\in L$, we have the following convergences: \begin{equation}\la{Jeudi14janv20160}\lim_{N \to \infty}\ff{N} \tr \bM_\ell\bM_{\ell'}-\ff{N} \tr \bM_\ell\cdot \ff{N} \tr\bM_{\ell'}=\eta_{\ell \ell'}\end{equation} and \begin{equation}\la{Jeudi14janv20161}\lim_{N \to \infty}\ff{N} \tr \bM_\ell\bM_{\ell'}^*-\ff{N} \tr \bM_\ell\cdot \ff{N} \tr\bM_{\ell'}^*=\bet_{\ell \ell'}.\end{equation}
\ent 
\end{assum}


Under this sole hypothesis, we first have the following result, focused on the case where the $\bM_\ell$'s all have null trace, \ie focused on the projections of the above  $\bA_k$'s onto the space of such matrices.

\begin{theo}\la{Th_Mainbis} Under Hypothesis \ref{Assum201611},  if, for each $\ell$, $\tr(\bM_{\ell})=0$, then the  finite-dimensional marginal distributions of
\bbe\la{17101516h0bis}
  \left(\tr \big( \bA_{k} \bM_{\ell}  \big)  \ri)_{k\in K, \ell\in L}
\ee 
converge to the ones of  a complex centered Gaussian vector $\big( \mathcal{G}_{k,\ell}\big)_{k\in K, \ell\in L}$ with covariance
$$
\E\left[\mathcal{G}_{k,\ell} {\mathcal{G}_{k',\ell'}}\right]  =   \eta_{\ell \ell'} \tau(k ,k')  ,\qquad\; \text{ and }\qquad\; \E\left[\mathcal{G}_{k,\ell} \ovl{\mathcal{G}_{k',\ell'}}\right]  =   \beta_{\ell \ell'} \tau(k ,\ovl{k'}).
$$
\end{theo}

\begin{rem}\label{remcent}Note that by  invariance of the distribution of $\bA$ under unitary conjugation, we have $$\E \tr \big( \bA_{k} \bM_{\ell}  \big)=\E\big(\ff{N}\Tr \bA_{k}\big)   \Tr\bM_\ell,$$ hence the random variables of \eqref{17101516h0bis} are centered and the ones of \eqref{17101516h0} below  rewrite $$ \tr \big( \bA_{k} \bM_{\ell}  \big)  -\E\big(\ff{N}\Tr \bA_{k}\big)   \Tr\bM_\ell .$$
\end{rem}

The following theorem gives the joint fluctuations of the projections of the $\bA_k$'s on null trace matrices   and of their traces.
\begin{assum}\la{Assum2016112}The  finite-dimensional marginal distributions of the process  $\big(\Tr \bA_k-\E\Tr\bA_k \big)_{k \in K}$ converge to those of a random centered vector $(\mathcal{T}_k)_{k\in K}$
and  for each $\ell\in L$, there is $\al_\ell\in \C$ \st  \begin{equation}\la{Jeudi14janv20162}\lim_{N \to \infty} \ff N \tr \bM_\ell  = \al_\ell.\end{equation}
\end{assum}

\begin{theo}\label{Th_Main}Under Hypotheses \ref{Assum201611} and \ref{Assum2016112}, the  finite-dimensional marginal distributions of \bbe\la{17101516h0}
  \left(\tr \big( \bA_{k} \bM_{\ell}  \big)  -\E \tr \big( \bA_{k} \bM_{\ell}  \big)\ri)_{k\in K, \ell\in L}
\ee 
converge to the ones of   $\big( \mathcal{G}_{k,\ell}+\al_{\ell}\mathcal{T}_{k}\big)_{k\in K, \ell\in L}$ where $\big( \mathcal{G}_{k,\ell}\big)_{k\in K, \ell\in L}$ is a complex centered Gaussian vector independent from $\big(\mathcal{T}_{k}\big)_{k\in K}$ and with covariance
$$
\E\left[\mathcal{G}_{k,\ell} {\mathcal{G}_{k',\ell'}}\right]  =   \eta_{\ell \ell'} \tau(k ,k')  ,\qquad\; \text{ and }\qquad\; 
\E\left[\mathcal{G}_{k,\ell} \ol{\mathcal{G}_{k',\ell'}}\right]  = \bet_{\ell \ell'} \tau ( k, \ovl{k'}).
$$
\end{theo}

A  direct consequence of this theorem is the asymptotic independence of the projections of the matrices $\bA_k$ onto the subspace of null trace matrices from  their projections onto the orthogonal of this subspace:
\begin{cor}\label{Th_Mainter115} Under Hypotheses \ref{Assum201611} and \ref{Assum2016112}, suppose that for any $\ell \in L$, we have $\Tr (\bM_\ell)=0$. Then the processes $$  \left(\tr \big( \bA_{k} \bM_{\ell}  \big)  \ri)_{k\in K, \ell\in L}\qquad\text{ and }\qquad \big(\Tr \bA_k-\E\Tr\bA_k \big)_{k \in K}$$ are asymptotically independent.
\end{cor}

\begin{rem}\label{no_second_order_freeness_supposed} 
It has been proved in \cite{Mingo2007} that unitary invariance implies second order freeness in many cases. However,   Theorems \ref{Th_Mainbis} and \ref{Th_Main}, as well as their corollaries,  cannot be deduced   from  the theory of second order freeness. The reason is that neither the random matrices $\bA_k$ nor the matrices $\bM_{\ell}$ are supposed to have a second order limit. Even in the case where the random matrices $\bA_k$ have a second order limit  (see Section 2.2), the ``test matrices" that  we consider (\ie the matrices $\bM_\ell$) are not supposed to have a second order limit, nor to be well approximated by matrices having a second order limit. For example, if  $$\bM_\ell=\sqrt{N}\times\text{(an elementary matrix)}$$ (a typical case of application of our results),   then for any $p\ge 2$, the sequence  $$\ff{N}\Tr(\bM_\ell\bM_\ell^*)^p=N^{p-1}$$ does not have any finite limit as $N\to\infty$, nor is bounded    (which would be required to prove our results as application of second order freeness). 
\end{rem}

\subsection{Second-order freeness implies fluctuations of matrix elements}\label{21021720h}As explained in Remark \ref{no_second_order_freeness_supposed}, our results do not follow from second order freeness theory. However, we shall see in the following corollary that  they   extend the consequences  of the existence of a second order limit   for  unitarily invariant matrix ensembles.  
 Let $\C\langle x_k,x_k^*,k\in K\rangle$ denote the space of polynomials in the non commutative variables $x_k,x_k^*$, indexed by $k\in K$. The following corollary follows directly from     Theorem \ref{Th_Main}.
\begin{cor}\label{sofifme} Let $ \big(\bA_k\big)_{k \in K}$ be a collection of    $N\ti N$ random  matrices which is   invariant by unitary conjugation and  which converges  in second order $*$-distribution to some family $a=  (a_k)_{k \in K}$ in $(\mathcal{A},\tau_1,\tau_2)$  as $N\to\infty$. Let $ \big(\bM_\ell\big)_{\ell \in L}$ be a collection non random   matrices satisfying \eqref{Jeudi14janv20160}, \eqref{Jeudi14janv20161} and \eqref{Jeudi14janv20162}.

Then the  finite-dimensional marginal distributions of 
\bbe
  \left(\tr \big( P(\bA) \bM_{\ell}  \big)  -\E \tr \big( P(\bA) \bM_{\ell}  \big)\ri)_{P\in\C\langle x_k,x_k^*,k\in K\rangle, \ell \in L}
\ee  converge to the ones of 
  a complex centered Gaussian vector $$\big( \mathcal{H}_{P,\ell}\big)_{P\in\C\langle x_k,x_k^*,k\in K\rangle, \ell \in L}$$ such that, for all $P, Q\in\C\langle x_k,x_k^*,k\in K\rangle$ and $\ell, \ell'\in L$,
\beq
\E{\mathcal{H}_{P,\ell} {\mathcal{H}_{Q, \ell'}}} & = & \al_{\ell} \al_{\ell'}\tau_2(P(a),Q(a)) +\eta_{\ell \ell'} \big(\tau_1(P(a) Q(a)) - \tau_1(P(a))\tau_1(Q(a)) \big),\\
\E{\mathcal{H}_{P,\ell} \ol{\mathcal{H}_{Q, \ell'}}} & = &\al_{\ell} \ol{\al_{\ell'}} \tau_2(P(a),Q(a)^*) + \bet_{\ell \ell'}\big(\tau_1(P(a) Q(a)^*) - \tau_1(P(a))\tau_1(Q(a)^*) \big) .\\
\eeq
\end{cor}
 \begin{rem}\la{list_applications_240217}The following matrices have been shown to converge in second order $*$-distribution:
  \bgt
   \ite[-] Wishart matrices and matrices of the type $\bU\bA\bV$ or $\bU\bA\bU^*$, with $\bU,\bV$ independent and Haar distributed  on $U(N)$  and $\bA$ deterministic with a limit spectral distribution \cite{mingo-nica04,mingo-speicher06, Mingo2007, mingo-piotr-collins-speicher07}
 \ite[-] GUE matrices or more generally   matrix models where the entries interact via a potential    \cite{johansson},
 \ite[-] Ginibre matrices \cite{Rider2006},
 \ite[-] random unitary matrices distributed according to  the Haar measure on the unitary group $U(N)$ \cite{Diaconis1994},
 \ite[-] matrices arising from the heat kernel measure on $U(N)$ \cite{Levy2010} and on $GL_N(\mathbb{C})$ \cite{ToddGuillaume}.
 \ent
  A consequence of Corollary~\ref{sofifme} is that any non commutative polynomial in independent random matrices taken from the list above has asymptotically Gaussian entries, which are independent modulo a possible symmetry.\end{rem}

\subsection{Left and right unitary invariant matrices}
Here is another corollary on random matrices invariant by left and right unitary multiplication.  
\begin{cor}\label{cor09040369875015}
Let $\bA  = \big(\bA_k\big)_{k \in K}$ be a collection of   $N\ti N$ random  matrices such that: 
\bgt
\item[(a')]  $\bA $ is  invariant, in law, by  left and right multiplication by unitary matrices: for    any non random unitary matrix $\bU$,  $\bA \stackrel{\text{law}}{=}\big(\bU\bA_k\big)_{k \in K}\stackrel{\text{law}}{=}\big(\bA_k \bU\big)_{k \in K}  $;
\item[(b')] for each $k$ and each $p,q$, $\ff N \tr |\bA_k|^{2p}$ is bounded in $L^q$ independently of $N$;
\item[(c')] for each $k,k'$, the sequence $\ff N \tr \bA_k \bA_{k'}^* $ converges in $L^2$ to some non random limits denoted $\tau(k,\ol{k'})$. 
\ent
Let $ \big(\bM_\ell\big)_{\ell \in L}$ be a collection non random   matrices satisfying \eqref{Jeudi14janv20160}, \eqref{Jeudi14janv20161} and \eqref{Jeudi14janv20162}. Then the  finite-dimensional marginal distributions of   
$$
\lf( \tr(\bA_{k} \bM_{\ell})\ri)_{k\in K, \ell\in L}
$$
converge to the ones of a complex centered Gaussian vector $\big( \mathcal{G}_{k,\ell}\big)_{k\in K, \ell\in L}$ with covariance 
\beq
\E \mathcal{G}_{k,\ell} \mathcal{G}_{k',\ell'} \ = \ 0, & \text{ and } &
\E \mathcal{G}_{k,\ell} \ol{\mathcal{G}_{k',\ell'}} \ = \ \bet_{\ell,\ell'}\tau(k,\ol{k'}).
\eeq
\end{cor}
The proof of this corollary is postponed to Section \ref{proof0905110520115}: we show that the hypotheses of the corollary imply Hypotheses \ref{Assum201611} and  \ref{Assum2016112}.

\subsection{Permutation matrix entries under randomized basis}
 Let $S$ be a uniform random $N\times N$ permutation matrix. For $T_d$ the number of $d$-cycles of the underlying permutation,  the   distribution of $(T_d)_{d\ge 1}$  converges as $N\to\infty$ to  a Poisson process $(\mathcal{Z}_d)_{d\ge 1}$ on the set of positive integers  with intensity $1/d$ (see \cite{Arratia1992}). It implies  that each trace $\Tr(S^k)$ ($k\ge 1$) converges in distribution to $\sum_{d|k}d\mathcal{Z}_d$. Thanks to Theorem \ref{Th_Main} and Remark~\ref{remcent}, we deduce directly  the following result about the matrix entries of a uniform permutation matrix $S$ conjugated by a uniform unitary matrix.

\begin{cor}
Let $S$ be a $N\times N$ random permutation matrix which is uniformly distributed, $U$ be a $N\times N$ random unitary matrix which is Haar distributed , and $ \big(\bM_\ell\big)_{\ell \in L}$ a collection of non random   matrices satisfying \eqref{Jeudi14janv20160}, \eqref{Jeudi14janv20161} and \eqref{Jeudi14janv20162}. Then the  finite-dimensional marginal distributions of   
$$
\lf( \tr\big(
US^kU^* \bM_{\ell}\big)\ri)_{k\ge 1, \ell\in L}
$$
converge to the ones of $\big( \mathcal{G}_{k,\ell}+\alpha_\ell \sum_{d|k}d\mathcal{Z}_d\big)_{k\ge 1, \ell\in L}$, where $\big( \mathcal{G}_{k,\ell}\big)_{k\ge 1, \ell\in L}$ is a complex centered Gaussian vector with covariance
\beq
\E \mathcal{G}_{k,\ell} \mathcal{G}_{k',\ell'} \ = \ 0, & \text{ and } &
\E \mathcal{G}_{k,\ell} \ol{\mathcal{G}_{k',\ell'}} \ = \ \one_{k=k'}\beta_{\ell,\ell'},
\eeq
and $(\mathcal{Z}_d)_{d\ge 1}$ is a Poisson process  on the set of positive integers  with intensity $1/d$ which is independent from  $\big( \mathcal{G}_{k,\ell}\big)_{k\in \mathbb{N}, \ell\in L}$.
\end{cor}

This is to be compared with the results of  \cite{Tsou2015}, where the  entries of the  matrix $S$ conjugated by a uniform random orthogonal matrix are studied. 

%

\subsection{Low rank perturbation for Gaussian elliptic matrices}

Matrices from the \emph{Gaussian elliptic ensemble}, first introduced in \cite{SCS}, can be defined as follows.
\begin{defi} A  \emph{Gaussian elliptic matrix of parameter $\rho \in [-1,1]$}  is a 
  random matrix  $\bY=[y_{ij}]_{i,j=1}^N$ \st \bgt
\ite $\{(y_{ij},y_{ji}), 1 \leq i < j \leq N\} \cup \{y_{ii}, 1 \leq i \leq N\}$ is a family of independent  random vectors,
\ite $\{(y_{ij},y_{ji}), 1 \leq i < j \leq N\}$ are i.i.d. Gaussian, centered,    such that
$$
\E y_{ij}^2 = \E y_{ji}^2 = \E y_{ij} \ol{y_{ji}} = 0, \quad \E |y_{ij}|^2 = \E |y_{ji}|^2 = 1 \ \quad \text{ and } \quad \E y_{ij}y_{ji} = \rho
$$
\ite $\{y_{ii}, 1 \leq i \leq N\}$ are i.i.d. Gaussian, centered,   such that
$$
\E y_{ii}^2 = \rho \quad \text{ and } \quad \E |y_{ii}|^2 = 1.
$$
\ent
\end{defi}

\begin{rem} \la{rem12581012201548}
Gaussian elliptic matrices can be seen as an interpolation between     \emph{GUE matrices} and   \emph{Ginibre matrices}. Indeed, a Gaussian elliptic matrix $\bY$ of parameter $\rho$   can   be realized  as 
$$
\bY \ = \ \sqrt{\f{1+\rho}{2}}\bH_1 + \ii \sqrt{\f{1-\rho}{2}}\bH_2,
$$
where $\bH_1$ and $\bH_2$ are two independent GUE matrices from the \emph{GUE}.  Hence     GUE matrices (resp.  Ginibre matrices) are    {Gaussian elliptic matrices} of parameter $1$ (resp. $0$).
\end{rem}

One can also define more general \emph{elliptic random matrices} (see \cite{Naum,NGUORO,or13,or14} for more details). In our case, it is easy to see (using for example Remark \ref{rem12581012201548}) that the {Gaussian elliptic ensemble} is invariant in distribution by unitary conjugation, which allows us to use our Theorem \ref{Th_Main} for this model. In this section, we are interested in     the outliers in the spectrum of these matrices. It is known (see \cite{SCS}) that when you renormalize  the matrix $\bY$ by $\sqrt N$, its limiting eigenvalue distribution is the uniform measure $\mu_\rho$ on the ellipse 
\beqy \la{ellipse2015481418sq}
\mathcal{ E}_\rho \ := \ \lf\{z \in \C \ ; \   \f{(\re z)^2}{(1+\rho)^2} + \f{(\im z)^2}{(1 - \rho)^2} \ \leq \ 1\ri\}.
\eeqy
Also, we know that adding  a finite rank matrix $\bP$ to such a matrix $\bY$ barely alters  its spectrum from the global point of view  (see \cite[Theorem 1.8]{NGUORO}),  but    may give rise to outliers. The generic location of the outliers has already been studied   (see \cite{or13}), but the authors did not consider the fluctuations. 

For all $N \geq 1$, let $\bX_N := \ff{\sqrt N} \bY_N$ where $\bY_N$ is an $N \ti N$  {Gaussian elliptic matrix} of parameter $\rho$ and let $\bP_N$ be a $N \ti N$ random matrix independent from $\bX_N$ whose rank is bounded by an integer $r$ (independent from $N$). We consider the additive pertubation
\beq
\wt\bX_N & := & \bX_N + \bP_N.
\eeq
Since, for any unitary matrix $\bU$ which is independent from $\bX_N$, we have $\bX_N \overset{(d)}{=} \bU \bX_N \bU^*$, we can assume that $\bP_N$ has the following block structure
\beq
\bP_N & = & \bpm \bP & 0 \\ 0 & 0 \epm , \ \ \ \ \ \text{ where } \bP \text{ is a } 2r \ti 2r \text{ matrix} 
\eeq
(indeed, any complex matrix is unitarily similar to a upper triangular matrix and since the rank of $\bP_N$ is lower than $r$, we have $\dim(\operatorname{Im} \bP_N + (\operatorname{Ker}\bP_N)^{\perp}) \leq 2r$). 

%
%

\begin{theo}[Outliers for finite rank perturbations of a Gaussian elliptic matrix]  \la{theo08521412020154}
Let $\ep>0$.  Suppose that $\bP_N$ does not have any eigenvalue $\lam$ such that
\beqy \la{0826141215896l2015} |\lam|\ > \ 1\quad\trm{ and }\quad
1 + |\rho| + \ep \ < \ \lf|\lam + \f{\rho}{\lam } \ri| \ < \ 1 + |\rho| + 3 \ep, 
\eeqy
and has exactly  $j \leq r$ eigenvalues $\lam_1(\bP_N),\ld,\lam_j(\bP_N)$ (counted with multiplicity) such that, for each $i=1, \ld, j$, 
\beqy \la{0826141215896l2015bis}  |\lam_i(\bP_N)|\ > \ 1\quad\trm{ and }\quad 
\lf|\lam_i(\bP_N) + \f{\rho}{\lam_i(\bP_N)} \ri| \ > \ 1 + |\rho| + 3 \ep.\eeqy
Then, with  probability tending to one, $\wt\bX_N := \bX_N + \bP_N$ possesses exactly $j$ eigenvalues $\wt\lam_1, \ld, \wt\lam_j$ in $\{z \in \C\,;\, |z|>1+|\rho| + 2\ep\}$ and after a proper labeling 
\begin{equation}\label{1411619h}
\wt\lam_i  \ = \ \lam_i(\bP_N) +  \f{\rho}{\lam_i(\bP_N)} + \oo1, 
\ee
for each $1 \leq i \leq j$.
\end{theo}
\begin{rem}
In \cite{or13}, the authors prove this result for \emph{real elliptic random matrices} and have a more precise statement. Indeed, they replace in our conditions \eqref{0826141215896l2015} and  \eqref{0826141215896l2015bis} the annulus  $\{z \in \C\,;\,  1 + |\rho| + \ep< |z|<1 + |\rho|+3\ep\}$ (resp.  $\{z \in \C, \ |z|> 1 + |\rho| + 3 \ep \}$) by $\mathcal{E}_{\rho,3\ep} \backslash \mathcal{E}_{\rho,\ep}$ (resp. by $\mathcal{E}_{\rho,3\ep}^c$) where $\mathcal{E}_{\rho,\ep}$ is a $\ep$-neighborhood of the ellipse $\mathcal{E}_\rho$ (see \eqref{ellipse2015481418sq}). Our proof relies on the identity
$$
\tr \inve{z-\bX } \ = \ \sum_{k \geq 0} z^{-k-1} \tr \bX^k,
$$ which is true only when $|z|$ is larger than the spectral radius of $X$, 
this is why \eqref{0826141215896l2015} and \eqref{0826141215896l2015bis} are circular conditions, instead of elliptic ones.
\end{rem}

To study the fluctuations of the outliers $\wt\lam_i$ around their generic locations as given by \eqre{1411619h}, we need to specify   the shape of the matrix $\bP$ as it is done in \cite{FloJean}. Indeed, since $\bP$ is not Hermitian, we need to introduce its Jordan Canonical Form (JCF) which is supposed to be independent  of $N$, but for its kernel part. We know that, in a proper basis, $\bP$ is a direct sum of \emph{Jordan blocks}, i.e. blocks of the form
\beqy
\bR_p(\tta) & = &  \bpm \tta & 1 & & (0) \\ 
                             &\ddots&\ddots &\\ 
                             \multicolumn{2}{c}{\multirow{2}{*}{(0)}} &\ddots & 1 \\
                             \multicolumn{2}{c}{}                     & &\tta \epm \ \in \ \C ^{p\ti p}, \qquad \tta \in \C, \, p \geq 1
\eeqy 

Let us denote by $\tta_1, \ld, \tta_q$ 
the distinct eigenvalues of $\bP $ satisfying condition \eqref{0826141215896l2015bis}. For convenience, we shall write from now on
\beqy
\hat \tta_i & := & \tta_i + \f{\rho}{\tta_i}.
\eeqy 
We introduce a positive integer $\al_i$, some positive integers      $p_{i,1}> \cdots> p_{i,\al_i}$ corresponding to the distinct sizes of the blocks relative to the eigenvalue $\tta_i$  and $\bet_{i,1}, \ld, \bet_{i, \al_i}$ \st  for all $j$, $\bR_{p_{i,j}}(\tta_i)$ appears $\bet_{i,j}$ times, so that, for a certain $2r \ti 2r$ invertible  matrix $\bQ$, we have:\\ 

\beqy\label{Eq141130042001587872015} \bJ &=& \bQ^{-1}\bP\bQ  \ = \
 \Bigg(  \bigoplus_{i=1}^q \ \bigoplus_{j=1}^{\alp_i} \! \underbrace{\begin{pmatrix}\bR_{p_{i,j}}(\tta_i)&&\\ &\ddots&\\ &&\bR_{p_{i,j}}(\tta_i)\end{pmatrix}}_{
 \bet_{{i,j}} \trm{ blocks}}\Bigg) \ \bigoplus \  \hat\bP
 \eeqy
where $\oplus$ is defined, for square block matrices, by $\mathbf{M}\oplus \mathbf{N}:=\bpm \mathbf{M}& 0\\ 0&\mathbf{N}\epm$ and $\hat\bP$ is a matrix whose  eigenvalues $\tta$ are such that $|\tta|<1$ or $ |\tta + \rho \tta^{-1}|<1 + \rho + \ep$.\\
\indent The asymptotic orders of the fluctuations of the eigenvalues of $\bX_N+\bP_N $ depend  on the sizes $p_{i,j}$ of the blocks. We know, by Theorem \ref{theo08521412020154}, that  there are $\sum_{j=1}^{\alp_i}p_{ij}\ti\bet_{i,j}$ eigenvalues ${\wt\lambda}$ of $\bX_N+\bP_N$  which tend  to $\hat{\tta_i} = \tta_i + \rho \tta_i^{-1}$: we shall write them with a $\hat{\tta_i}$ 
 on the top left corner, as follows $$\lexp{\hat\tta_{i}}{{\wt\lambda}}.$$ Theorem \ref{thbelow0911141210002015} below will state that,  for each block with size 
 $p_{i,j}$  corresponding to $\tta_i$  in  the JCF of $\bP$, there are $p_{i,j}$ eigenvalues (we shall write them with $p_{i,j}$ on the bottom left corner: $\lbinom{\hat\tta_{i}}
 {p_{i,j}}{{\wt\lambda}}$) whose  convergence rate will be $N^{-1/(2p_{i,j})}$. As there are $\bet_{{i,j}}$ blocks of size 
 $p_{i,j}$, there are actually $p_{i,j}\times \bet_{{i,j}}$ eigenvalues tending to $\hat\tta_i$ with  convergence rate    
 $N^{-1/(2p_{i,j})}$ (we shall write them $\lbinom{\hat\tta_i}{p_{i,j}}{{\wt\lambda}_{s,t}}$ with $s \in \{1,\ldots,p_{i,j}\}$ and 
 $t \in \{1,\ldots,\bet_{{i,j}}\}$). It would be convenient to denote by $\Lambda_{i,j}$ the vector with size $p_{i,j}\times \bet_{{i,j}}$ defined by 
 \beqy \label{defLambda300420152015}
 \Lambda_{i,j} \egd  \displaystyle \left(N^{1/(2p_{i,j})} \cdot \Big(\lbinom{\hat\tta_{i}}{p_{i,j}}{ {\wt\lambda}_{s,t}} - \hat\tta_i \Big) \right)_{\substack{1\le s\le p_{i,j}\\ 1\le t\le \bet_{i,j}}}.
\eeqy

\indent As in \cite{FloJean}, we define now the family of random matrices that we shall use to characterize the limit distribution of the $\Lambda_{i,j}$'s. 
For each $i=1, \ldots, q$, let $I(\tta_i)$ (resp. $J(\tta_i)$) denote the set, with cardinality $\sum_{j=1}^{\al_i}\bet_{i,j}$, of indices in $\{1, \ld,2 r\}$ corresponding to the first (resp. last) columns of the blocks $\bR_{p_{i,j}}(\tta_i)$ ($1\le j\le \al_i$) in \eqre{Eq141130042001587872015}.
\begin{rem} \la{rem16312909201420152015}
 Note that the columns of $\bQ$ (resp. of $(\bQ^{-1})^*$) whose index belongs to $I(\tta_i)$ (resp. $J(\tta_i)$) are   eigenvectors of $\bP$  (resp. of $\bP^*$) associated to $\tta_i$ (resp. $\ol{\tta_i}$). See \cite[Remark 2.8]{FloJean}.
 \end{rem}

Now, let  \bbe\la{2121414h220152015}
\lf({m}^{\tta_i}_{k,\ell}\ri)_{\substack{1 \leq i \leq q_{\hphantom{1.}} \quad\quad \ \ \\   (k,\ell)\in J(\tta_i)\ti I(\tta_i)}} \ee 
 be the random centered complex Gaussian vector with covariance 
 
\beqy \la{100714211412015}
\begin{array}{rcl} 
\Ec{m_{k,\ell}^{\tta_i} m_{k',\ell'}^{\tta_{i'}}} & = &\ds \Big(\ff{\tta_i \tta_{i'}-\rho} - \ff{\tta_i \tta_{i'}} \Big)\delta_{k,\ell'}\delta_{k',\ell} \\ 
 \Ec{m_{k,\ell}^{\tta_i} \ol{ m_{k',\ell'}^{\tta_{i'}}}} & = &\ds \Phi_\rho(\hat{\tta_i},\hat{\tta}_{{i'}})\be_k \bQ^{-1} (\bQ^{-1})^* \be_{k'} \cdot \be_{\ell'} \bQ^* \bQ \be_{\ell} ,
 \end{array}
\eeqy 
where    $\be_1, \ld, \be_{2r}$ are  the column vectors of the canonical basis of $\C^{2r}$ and 
\beq
\Phi_\rho(z,z') & = & \int \ff{z-w}\ff{\ol{ z'} - \ol w} \mu_\rho (\mathrm{d}w) - \int \ff{z-w} \mu_\rho (\mathrm{d}w)\int \ff{\ol{ z'} - \ol w} \mu_\rho (\mathrm{d}w).
\eeq

\begin{rem}
In Section \ref{section141212015fga}, the random vector of  \eqref{2121414h220152015} will appear as limit in the   convergence:     \beq
 \lf(\sqrt{N}\be_k^* \bQ^{-1}\lf( \inve{\hat\tta_i-\bX_N} - \ff {\tta_i}\ri)\bQ\be_\ell\ri)_{\substack{1 \leq i \leq q_{\hphantom{1.}} \quad\quad \ \ \\  (k,\ell)\in J(\tta_i)\ti I(\tta_i)}} \ {\cloi} \ \lf({m}^{\tta_i}_{k,\ell}\ri)_{\substack{1 \leq i \leq q_{\hphantom{1.}} \quad\quad \ \ \\  (k,\ell)\in J(\tta_i)\ti I(\tta_i)}}\,.
    \eeq 
 This convergence is a consequence of Theorem \ref{Th_Main}.\\
\end{rem}

\begin{rem}
When $\rho = 0$, one has
\beq
\Phi_0(z,z') & = & \ff{\pi} \int_{|w|\leq 1} \ff{z-w}\ff{\ol{ z'} - \ol w} \mathrm{d}w - \ff{\pi} \int_{|w|\leq 1} \ff{z-w} \mathrm{d}w \ff{\pi} \int_{|w|\leq 1} \ff{\ol{ z'} - \ol w} \mathrm{d}w  \\ &= & \ff{z\ol{z'} -1} - \ff{z \ol{z'}} \ = \ \ff{z\ol{z'}(z \ol{z'}-1)}. 
\eeq
We recover the expression of the covariance in the Ginibre case (see \cite{FloJean}). Also, the expression of $\Phi_1$ corresponds to the covariance in the GUE case (see \cite{Rochet}).
\end{rem}

 For each $i,j$, let $K(i,j)$ (resp. $K(i,j)^-$) be the set, with cardinality $\bet_{i,j}$ (resp. $\sum_{j'=1}^{j-1}\bet_{i,j'}$), of  indices in $J(\tta_i)$    corresponding to a block of the type $\bR_{p_{i,j}}(\tta_i)$ (resp. to a block of the type  $\bR_{p_{i,j'}}(\tta_i)$ for  $j'<j$). In the same way, let $L(i,j)$ (resp. $L(i,j)^-$) be the set, with the same cardinality as $K(i,j)$ (resp. as $K(i,j)^-$), of indices in $I(\tta_i)$ corresponding to a block of the type $\bR_{p_{i,j}}(\tta_i)$ (resp. to a block of the type  $\bR_{p_{i,j'}}(\tta_i)$ for  $j'<j$).  Note that $K(i,j)^-$ and $L(i,j)^-$ are empty if $j=1$.  Let us define the random matrices 
 \beqy \label{section2.320152015}
 \ope{M}^{\tta_i,\mathrm{I}}_{j}\egd[m^{\tta_i,n}_{k,\ell}]_{\ds^{k\in K(i,j)^-}_{\ell \in L(i,j)^-}} &\qquad \qquad \qquad& \ope{M}^{\tta_i,\dI\dI}_{j}\egd[m^{\tta_i}_{k,\ell}]_{\ds^{k\in K(i,j)^-}_{\ell\in L(i,j)}} \nonumber\\
 &&\\
 \ope{M}^{\tta_i,\dI\dI\dI}_{j,n}\egd[m^{\tta_i}_{k,\ell}]_{\ds^{k\in K(i,j)}_{\ell\in L(i,j)^-}} &\qquad \qquad \qquad& \ope{M}^{\tta_i,\dI\dV}_{j}\egd[m^{\tta_i,n}_{k,\ell}]_{\ds^{k\in K(i,j)}_{\ell \in L(i,j)}} \nonumber
 \eeqy
 and then let us define the  matrix 
${\bM}^{\tta_i}_{j}$   as 
\bbe\la{21021415h2015}\bM^{\tta_i}_{j}\egd  \tta_i \left(\ope{M}^{\tta_i,\dI\dV}_{j}-\ope{M}^{\tta_i,\dI\dI\dI}_{j}\inve{\ope{M}^{\tta_i,\dI}_{j}}\ope{M}^{\tta_i,\dI\dI}_{j} \right)
\ee
\begin{rem} \label{rem23021420152015}
It follows from the fact that the matrix $\bQ$ is invertible  that $\ope{M}^{\tta_i,\dI}_{j}$ is a.s. invertible and that so is $\bM^{\tta_i}_{j}$. 
\end{rem}

 \indent Now, we can state the result about the fluctuations of the outliers.

\begin{theo} \la{thbelow0911141210002015}
\begin{enumerate}
\item As $N$ goes to infinity, the random vector  $$\displaystyle \left(\Lambda_{i,j} \right)_{\substack{1 \leq i \leq q_{\hphantom{1}} \\ 1 \leq j\leq \alp_i } } $$ defined at \eqre{defLambda300420152015} converges  in distribution to  a  random vector $$\displaystyle \left(\Lambda^\infty_{i,j} \right)_{\substack{1 \leq i \leq q_{\hphantom{1}} \\ 1 \leq j\leq \alp_i } } $$ with joint distribution defined by the fact that,
for each $1 \leq i \leq q$ and $1 \leq j \leq \al_i$,  $\Lambda_{i,j}^\infty$ is the collection of the  ${p_{i,j}}^{\trm{th}}$ roots of the  eigenvalues of the random matrix $\bM^{\tta_i}_{j}$. 
\item The distributions of the  random matrices $\bM^{\tta_i}_{j}$ are absolutely continuous with respect to the Lebesgue measure  and  the  random vector $\displaystyle \left(\Lambda^\infty_{i,j} \right)_{\displaystyle^{1 \leq i \leq q}_{1 \leq j\leq \alp_i}} $   has  no deterministic coordinate. 
\end{enumerate}
\end{theo}

\begin{rem} Each non zero complex number has exactly $p_{i,j}$ ${p_{i,j}}^{\trm{th}}$ roots, drawing a regular $p_{i,j}$-sided polygon. Moreover, by the second part of the theorem, the spectrums of the $\bM^{\tta_i}_{j}$'s almost surely do not contain $0$,  so  each $\Lambda_{i,j}^\infty$ is actually a complex random vector with $p_{i,j}\ti \bet_{i,j}$ coordinates, which draw $\bet_{i,j}$ regular $p_{i,j}$-sided polygons.

\end{rem}

\begin{rem} \la{rem1519kl2901}
We notice that in the particular case where the matrix $\bQ$ is unitary, the covariance of the Gaussian variables $\lf({m}^{\tta_i}_{k,\ell}\ri)_{\substack{1 \leq i \leq q_{\hphantom{1.}} \quad\quad \ \ \\   (k,\ell)\in J(\tta_i)\ti I(\tta_i)}}$ can be rewritten 
\beq
\Ec{m_{k,\ell}^{\tta_i} m_{k',\ell'}^{\tta_{i'}}} & = &\ds \Big(\ff{\tta_i \tta_{i'}-\rho} - \ff{\tta_i \tta_{i'}} \Big)\delta_{k,\ell'}\delta_{k',\ell}, \\ 
 \Ec{m_{k,\ell}^{\tta_i} \ol{ m_{k',\ell'}^{\tta_{i'}}}} & = &\ds \Phi(\hat{\tta}_i, \hat{\tta}_{i'})\delta_{k,k'} \delta_{\ell,\ell'}.
\eeq 
Which means that for any $i,i'$ such that $i \neq i'$, the familly $\lf({m}^{\tta_i}_{k,\ell}\ri)_{   (k,\ell)\in J(\tta_i)\ti I(\tta_i)}$ is independent from $\lf({m}^{\tta_{i'}}_{k,\ell}\ri)_{   (k,\ell)\in J(\tta_{i'})\ti I(\tta_{i'})}$. Indeed, since the Jordan blocks associated to $\tta_i$ are distinct from  those associated to $\tta_{i'}$,  the sets $I(\tta_i)$ and $J(\tta_i)$ don't share any common index with  $I(\tta_{i'})$ and $J(\tta_{i'})$. We can deduce that in this particular case, all the fluctuations around $\tta_i$ are independent from those around $\tta_{i'}$ (see \cite[section 2.3.1.]{FloJean} for more details). \\
However, in the general case, there is no particular reason to have independance between the fluctuations around two spikes at macroscopique distance. To illustrate this phenomenon, we can take the same particular example than \cite[Example 2.17]{FloJean} since a Ginibre matrix is also a Gaussian elliptic matrix. In this example, the authors of \cite{FloJean} took a matrix $\bP$ of the form
$$
\bP \ = \ \bQ\bpm \tta & 0 \\ 0 & \tta' \epm \bQ^{-1}, \ \quad \quad \ \bQ \ = \ \bpm 1 & \kappa \\ \kappa & 1 \epm, \ \ \ \kappa \neq \pm 1,
$$
and they empirically confirmed that, in the case $\kappa \neq 0$, the fluctuations of the outliers around $\tta$ are correlated with these around $\tta'$.


\end{rem}

\section{Proofs of Theorem \ref{Th_Mainbis} and Theorem \ref{Th_Main}}\label{proofs}

\subsection{Preliminary result}

\noindent Let $  \big(\bB_k\big)_{k \in K}$ be a collection of (implicitly depending on $N$)  $N\ti N$ random  matrices such that
\bgt
\item[(i)] for each $k\in K$, almost surely,  $\tr \bB_k=0$;
\item[(ii)] for each $k\in K$, and each $p,q\geq 1$,  $\ff{N} \tr |\bB_k|^{2p}$ is bounded in $L^q$ independently of $N$;
\item[(iii)] for each $k,k'\in K$, we have the following convergences to nonrandom variables in $L^2$
$$\lim_{N \to \infty}\ff{N} \tr \bB_k\bB_{k'}=\tau(k,k')\ \text{ and }\lim_{N \to \infty}\ff{N} \tr \bB_k\bB_{k'}^*=\tau(k,\ovl{k'}).$$
\ent
 Let also $  \big(\bM_\ell\big)_{\ell \in L}$ be a collection non random   matrices such that
 \bgt
\item[(iv)] for each $\ell,\ell'\in L$, we have the following convergences $$\lim_{N \to \infty}\ff{N} \tr \bM_\ell\bM_{\ell'}-\ff{N} \tr \bM_\ell\cdot \ff{N} \tr\bM_{\ell'}=\eta_{\ell \ell'}$$ and $$\lim_{N \to \infty}\ff{N} \tr \bM_\ell\bM_{\ell'}^*-\ff{N} \tr \bM_\ell\cdot \ff{N} \tr\bM_{\ell'}^*=\bet_{\ell \ell'}.$$
\ent
 At last,   let $\bU=\bU^{(N)}$ be an $N\ti N$ Haar-distributed unitary random matrix independent of $\big(\bB_k\big)_{k \in K}$.

\begin{propo}\la{th_mainprelim}
Let us fix $p\ge 1$,  $(k_1,\ldots,k_p) \in K^p$ and $(\ell_1,\ldots,\ell_p)\in L^p$.
If $(i)$, $(ii)$, $(iii)$ and $(iv)$ hold, then the centered vector
\bbe\la{17101516h}
  \left(\tr \big(\bU\bB_{k_i}\bU^{*}\bM_{\ell_i}  \big)  \ri)_{1 \leq i\leq p}
\ee 
converges in distribution, as $N\to\infty$,  to a complex centered Gaussian vector $\big( \mathcal{G}_{i}\big)_{1 \leq i \leq p}$ such that, for all $i,i'$,
$$
\E{\mathcal{G}_{i} {\mathcal{G}_{i'}}} =  \eta_{\ell_i \ell_{i'}} \tau({k_i},{k_{i'}}) ,\ \text{ and }\
\E{\mathcal{G}_{i} \ol{\mathcal{G}_{i'}}}  = \bet_{\ell_i \ell_{i'}}  \tau({k_i},\ovl{k_{i'}}).
$$
Besides, for any sequence $(Y_N)$ of bounded random variables such that \bgt\ite $Y_N$ is independent of $\bU^{(N)}$, \ite  $\E Y_N$ has a   limit $L_Y$\ent  and any polynomial $f$ in $p$ complex   variables and their conjugates, we have  $$\lim_{N\to\infty} \E \lf[Y_Nf\lf( \tr \big(\bU\bB_{k_i}\bU^{*}\bM_{\ell_i}  \big) ,1 \leq i\leq p \ri)\ri]\;=\;L_Y\E \lf[f\lf(    \mathcal{G}_{i}  ,1 \leq i\leq p \ri)\ri] .$$
\end{propo}



\bpr First, we can suppose the $\bB_k$'s and the $\bM_\ell$'s are all Hermitian (which makes the entries of the vector of \eqre{17101516h} real), up to changing 
\beq\lf( \bB_k \ri)_{k \in K} & \longrightarrow & \lf( \bB_{(k,1)}:=\ff{2}\big(\bB_k+\bB_k^{*}\big) ,\;\; \bB_{(k,2)}:= \ff{2\mathrm{i}} \big(\bB_k  -\bB^{*}_k \big) \ri)_{(k,\ep) \in K\ti\{1,2\}},\\
\lf( \bM_\ell \ri)_{\ell \in L} & \longrightarrow & \lf( \bM_{(\ell,1)}:=\ff{2}\big(\bM_\ell+\bM_\ell^{*}\big) , \;\;\bM_{(\ell,2)}:= \ff{2\mathrm{i}} \big(\bM_\ell  -\bM^{*}_\ell \big) \ri)_{(\ell,\ep) \in L\ti\{1,2\}}.
\eeq
 Second, as all $\bB_k$'s have null trace, up to changing $\bM_\ell\to\bM_\ell-\ff{N}\tr \bM_\ell$, one can suppose that all $\bM_\ell$'s have null trace.

 To prove the full proposition, it suffices to prove the convergence 
 $$\lim_{N\to\infty} \E \lf[Y_N\prod_{i=1}^n \tr \big(\bU\bB_{k_i}\bU^{*}\bM_{\ell_i}  \big) \ri]\;=\;L_Y\E \lf[\prod_{i=1}^n \mathcal{H}_{i} \ri] .$$
 for any $n\ge 1$,  $(k_1,\ldots,k_n) \in K^n$, $(\ell_1,\ldots,\ell_n)\in L^n$ and any sequence $(Y_N)$ of bounded random variables independent from $\bU^{(N)}$ such that $\lim_{N\to\infty}\E Y_N=L_Y$. Indeed, we can take each $k$ as many times as we want in $(k_1,\ldots,k_n)$ (and the same for $\ell$), which implies the convergence of the expectation of any polynomials as wanted and consequently the convergence in distribution of finite dimensional marginals.
 
Let $n \ge 1$, and $\perm_n$ be the $n$-th symmetric group, and $\perm_{n,2}$ be the subset of permutations in $\perm_n$ with only cycles of length $2$. We denote be $\#\sigma$ the number of cycles of $\sigma\in \perm_n$ and by $\fix(\sigma)$ the number of fixed points of $\si$. The neutral element of $\perm_n$ is denoted by $\id_n$.
For any $\si \in \perm_n$, we set 
$$
\tr_\si \big( \bN_i\big)_{i=1}^n \ = \ \prod_{\substack{(t_1 t_2 \cdots t_m)\\\text{cycle of }\sigma}}\tr \big( \bN_{t_1} \bN_{t_2}\cdots \bN_{t_m}\big)
$$
For example, for $\si \in \perm_6$, $\si := (1,2,3,4,5,6) \ \mapsto \ (3,2,4,1,6,5)$
$$
\tr_\si \big( \bN_i\big)_{i=1}^6 \ = \ \tr\big( \bN_1 \bN_3 \bN_4\big) \tr \big( \bN_2\big) \tr \big( \bN_5 \bN_6\big).
$$

\begin{lemme}\la{1810201511}
Let $n\ge 1$, $(k_1,\ldots,k_n) \in K^n$, $(\ell_1,\ldots,\ell_n)\in L^n$, and $(Y_N)$ be any sequence of bounded random variables such that $\lim_{N\to\infty}\E Y_N=L_Y$. With the above assumptions on $(\bM_\ell)_{\ell\in L}$ and $(\bB_k)_{k\in K}$, we have, for all $\gamma$ and $\si$ in $\perm_n$,$$
\tr_\gamma \big( \bM_{\ell_{i}}\big)_{i=1}^n \ = \  \one_{\fix(\ga)=0}\, \times\,\OO{N^{n  /2}}
$$
and
$$
\E\left[ Y_N{\tr_\si \big( \bB_{k_{i}} \big)_{i=1}^n}\right] \ = \ \one_{\sigma\in\perm_{n,2}}N^{n/2}L_Y \prod_{\substack{(i,j)\\\text{cycle of }\sigma}}\tau(k_i,k_j) +o(N^{n/2}).
$$
\end{lemme}
\bpr Because $\bB_k$'s and   $\bM_\ell$'s   have null traces, the formulas are true in presence of fixed points. Thus, we can assume that $\sigma$ and $\gamma$ have no fixed point.

The first result comes from Lemma \ref{lem101206102015} and from the fact that, for each $\ell$, 
$
  \tr \bM_{\ell }^2 \ = \ \OO{N}.  
$

The second result can be proved in two steps. First, if $\sigma\notin\perm_{n,2}$, the non-commutative H\"older's inequality (see \cite[Appendix A.3]{agz}) and Hypothesis (ii) say us that
$$\left|\E\left[ Y_N{\tr_\si \big( \bB_{k_{i_j}} \big)_{j=1}^n}\right]\right| =O(N^{\#\sigma})=o(N^{n/2}).$$
If $\sigma\in\perm_{n,2}$ (and $n>0$ is even), we decompose $\sigma$ in $2$-cycles $
\si \ = \  (i_1 \ j_1) \cdots (i_{n/2} \ j_{n/2})$. By classical H\"older's inequality, the absolute difference between
$$N^{-n/2}\E\left[ Y_N{\tr_\si \big( \bB_{k_{j}} \big)_{j=1}^n}\right]=\E\left[ Y_N\prod_{t=1}^{n/2}\frac{1}{N} \tr \bB_{k_{i_t}}\bB_{k_{j_t}}\right] \text{and } \E\left[ Y_N\prod_{t=1}^{n/2-1}\frac{1}{N} \tr \bB_{k_{i_t}}\bB_{k_{j_t}}\right]\E\left[\frac{1}{N}\tr \bB_{k_{i_1}}\bB_{k_{j_1}}\right]$$
is less than
$$\E[Y_N^{2n}]^{1/2n}\prod_{t=1}^{n/2-1}\E\left[\left(\frac{1}{N} \tr \bB_{k_{i_t}}\bB_{k_{j_t}}\right)^{2n}\right]^{1/2n}\var\left(\frac{1}{N}\tr \bB_{k_{i_1}}\bB_{k_{j_1}}\right)$$
and consequently converges to $0$ - using again the non-commutative H\"older's inequality (see \cite[Appendix A.3]{agz}) and Hypothesis (ii) to control $\E(\frac{1}{N} \tr \bB_{k_{i_t}}\bB_{k_{j_t}})^{2n}$.
By a direct induction on $n/2$, it means that the expectation of product
$$N^{-n/2}\E\left[ Y_N{\tr_\si \big( \bB_{k_{j}} \big)_{j=1}^n}\right]=\E\left[ Y_N\prod_{t=1}^{n/2}\frac{1}{N} \tr \bB_{k_{i_t}}\bB_{k_{j_t}}\right]$$
has the same limit as the product of expectation
$$\E\left[ Y_N\right]\prod_{t=1}^{n/2}\E\left[\frac{1}{N} \tr \bB_{k_{i_t}}\bB_{k_{j_t}}\right],$$
and the result follows.
\epr
Let $n\ge 1$, $(k_1,\ldots,k_n) \in K^n$, $(\ell_1,\ldots,\ell_n)\in L^n$, and $(Y_N)$ be any sequence of bounded random variables such that $\lim_{N\to\infty}\E Y_N=L$.
Using \cite[Proposition 3.4]{Mingo2007} (and, first, an integration with respect to the randomness of $\bU$, and then a ``full expectation"),   we have
\beqy \label{11090611102015}
\E \lf[Y_N\prod_{i=1}^n \tr \big(\bU\bB_{k_i}\bU^{*}\bM_{\ell_i}  \big) \ri] & = & \sum_{\si,\gamma \in \perm_n} \wg(\si \gamma^{-1}) 
\E [Y_N\tr_\si \big( \bB_{k_i}\big)_{i=1}^n ]\tr_\gamma \big( \bM_{\ell_i}\big)_{i=1}^n,
\eeqy
where $\wg$ is the Weingarten function. We know from \cite[Coro. 2.7]{COL} and  \cite[Propo 23.11]{ns06} that, for any $\tau \in \perm_n$,
\begin{equation*}\la{eqWg171015}
\wg(\tau) \ = \ \OO{N^{\#\tau-2n}}\text{ and } \wg(1_n) \ = \ N^{-n}+O(N^{-n-2}).
\end{equation*}It implies, by Lemma \ref{1810201511},  that for $\si,\gamma \in \perm_n$,
 $$
\wg(\si \gamma^{-1}) 
\E [Y_N\tr_\si \big( \bB_{k_i}\big)_{i=1}^n ]\tr_\gamma \big( \bM_{\ell_i}\big)_{i=1}^n \ = \  \one_{\si\in\perm_{n,2}}\one_{\fix(\ga)=0}\,  O(N^{(\#(\si \gamma^{-1})-n)}   )=\one_{\gamma=\si\in\perm_{n,2}}  O(1),
$$
and more precisely, using the exact asymptotic for $\gamma=\si\in\perm_{n,2}$, that 
$$
\wg(\si \gamma^{-1}) 
\E [Y_N\tr_\si \big( \bB_{k_i}\big)_{i=1}^n ]\tr_\gamma \big( \bM_{\ell_i}\big)_{i=1}^n \ =\one_{\gamma=\si\in\perm_{n,2}} \ L_Y  \prod_{\substack{(i,j)\\\text{cycle of }\sigma}}\tau(k_i,k_j)\eta_{\ell_i\ell_j}+o(1).$$
As a consequence, we can rewrite \eqre{11090611102015} as
$$ \label{11090611102015prime}
\E \lf[Y_N\prod_{i=1}^n \tr \big(\bU\bB_{k_i}\bU^{*}\bM_{\ell_i}  \big) \ri]  = L_Y \sum_{ \si \in  {\perm}_{n,2}} \prod_{\substack{(i,j)\\\text{cycle of }\sigma}}\tau(k_i,k_j)\eta_{\ell_i\ell_j}+o(1),
$$
which is the wanted convergence in order to prove the proposition, since
$$\E \lf[\prod_{i=1}^n \mathcal{G}_{i} \ri]=\sum_{ \si \in  {\perm}_{n,2}} \prod_{\substack{(i,j)\\\text{cycle of }\sigma}}\tau(k_i,k_j)\eta_{\ell_i\ell_j}.$$
\epr

\subsection{Proof of Theorem \ref{Th_Main}}
First, note that, for each $k\in K$, $\E\bA_{k}=\ff{N}\E[\Tr(\bA_k)]I$, hence for $\bB_k:=\bA_k-\ff{N}\Tr(\bA_k)I$ and $T_k:=\ff{N}\Tr(\bA_k)-\E\ff{N}\Tr(\bA_k)$, one can write  $$\big(\bA_k-\E\bA_k\big)_{b\in K}= \big(\bB_k+T_kI\big)_{b\in K}.$$Let us now introduce an (implicitly depending on $N$) Haar-distributed unitary matrix $\bU$ independent of the collection $\bA$. By unitary invariance, we get $$\big(\bA_k-\E\bA_k\big)_{b\in K}= \big(\bB_k+T_kI\big)_{b\in K}\stackrel{\text{law}}{=} \big(\bU\bB_k\bU^*+T_kI\big)_{b\in K}.$$
Then, by  Proposition \ref{th_mainprelim}, we know that, for any $n\ge 1$, any $k_1, \ld, k_n\in K$ and any $\ell_1, \ld, \ell_n\in L$, the random vector $$\big(\Tr(\bU\bB_{k_i}\bU^*\bM_{\ell_i})\big)_{1\le i\le n}$$ converges in distribution to a complex centered Gaussian vector $\big( \mathcal{H}_{i}\big)_{1 \leq i \leq n}$ such that, for all $i,i'$,
\beq
\E{\mathcal{H}_{i} {\mathcal{H}_{i'}}} & = & \big(\tau({k_i} ,{k_{i'}}) - \tau ( {k_i})\tau ( {k_{i'}}) \big) \big(\eta_{\ell_i \ell_{i'}}-\al_{\ell_i} \al_{\ell_{i'}}  \big) ,  \\
\E{\mathcal{H}_{i} \ol{\mathcal{H}_{i'}}} & = &\big(\tau({k_i} ,\ovl{k_{i'}}) - \tau ( {k_i})\ovl{\tau ( {k_{i'}})} \big) \big(\eta_{\ell_i \ell_{i'}}-\al_{\ell_i} \ovl{\al_{\ell_{i'}} } \big). \eeq
Besides,  Proposition \ref{th_mainprelim} also says that $$\big(\Tr(\bU\bB_{k_i}\bU^*\bM_{\ell_i})\big)_{1\le i\le n}$$  is asymptotically independent from $\big(T_{k_i}\Tr(\bM_{\ell_i})\big)_{1\le i\le n}$, which converges in distribution, by Hypothesis (d), to $\big(\al_{\ell_i}\mathcal{T}_{k_i} \big)_{1\le i\le n}$. As it is clear, from the covariance of  $\big( \mathcal{G}_{i}\big)_{1 \leq i \leq n}$, that  for  $\big( \mathcal{H}_{i}\big)_{1 \leq i \leq n}$ independent from $\big(\al_{\ell_i}\mathcal{T}_{k_i}\al_{\ell_i}\big)_{1\le i\le n}$, we have $$\big( \mathcal{G}_{i}\big)_{1 \leq i \leq n}  \stackrel{\text{law}}{=} \big( \mathcal{H}_{i}\big)_{1 \leq i \leq n}+\big(\al_{\ell_i}\mathcal{T}_{k_i}\big)_{1\le i\le n},$$
the theorem is proved.\hfill$\square$

\subsection{Proof of Theorem \ref{Th_Mainbis}}
It is a direct application of Proposition \ref{th_mainprelim} since if $\tr \bM = 0$, then
$$
\tr (\bA \bM) \ = \ \tr \Big[\big( \bA - \ff N( \tr \bA )\bI\big)\bM\Big],
$$
so that one can assume that $\tr \bA = 0$.\hfill$\square$

\subsection{Proof of Corollary \ref{cor09040369875015}} \label{proof0905110520115}
We just need to show that the hypotheses of the corollary imply Hypotheses \ref{Assum201611} and  \ref{Assum2016112}.

The proof of Hypothesis \ref{Assum201611} comes down to    the following computations, where 
we introduce a Haar-distributed unitary matrix $\bU$ independent of $(\bA_k)_{k\in K}$ and use Equation (33) of  \cite{FloJean2}.
We have 
 \beq
 \E\big|\ff N\tr \bA_k \big|^2 & = & \ff {N^2} \Ec{ \E_\bU \big[ \tr (\bU \bA_k) \tr (\bA_k^*\bU^* )\big]} \\
& =& \ff {N^3}(1+\oo1) \Ec{\tr (\bA_k \bA_k^*)}  \ = \ \OO{\ff {N^2}} ,
 \eeq
and
\beq
\E{\big|\ff N \tr(\bA_k \bA_{k'})\big|^2} & =& \ff{N^2} \Ec{\E_\bU\big[\tr(\bA_k\bU\bA_{k'}\bU) \tr(\bA_{k'}\bU^*\bA_{k}\bU^*)\big]} \\
& = & \ff{N^4}(1+\oo1) \E\big( \tr(\bA_k \bA_k^*) \tr(\bA_{k'}\bA_{k'}^* )+ \tr(\bA_k \bA_{k'}^*) \tr(\bA_{k'}\bA_{k}^*)\big) \\
& = & \OO{\ff{N^2}}.
\eeq
Now, in order to show Hypothesis \ref{Assum2016112},  we want to prove that, for any fixed $r$, $\big(\tr(\bA_{k_1}),\ld,\tr(\bA_{k_r}) \big)_{i=1}^r$ is asymptotically Gaussian. Let $n\ge 1$ and $i_1,j_1,\ld,i_n,j_n\in \{1,\ld,r\}$, using \cite[Proposition 3.4]{Mingo2007}, we have

\beq
 \Ec{\prod_{\ell=1}^n \tr(\bA_{k_{i_\ell}}) \tr(\bA_{k_{j_\ell}}^*)} & = &  \Ec{\E_\bU\prod_{\ell=1}^n \tr(\bU\bA_{k_{i_\ell}}) \tr(\bA_{k_{j_\ell}}^*\bU^*)} \\
& =& \ff{N^n} \sum_{\si \in\perm_n} \Ec{\prod_{\ell=1}^n {\tr\big(\bA_{k_{i_\ell}} \bA_{k_{j_{\si(\ell)}}}^*\big)}}+\oo1.
\eeq
Then, one can prove that 
\beqy \label{115420112015}
\ff{N^n} \sum_{\si \in\perm_n} \Ec{\prod_{\ell=1}^n {\tr\big(\bA_{k_{i_\ell}} \bA_{k_{j_{\si(\ell)}}}^*\big)}} \ = \ \ff{N^n} \sum_{\si \in\perm_n} \prod_{\ell=1}^n \Ec{\tr\big(\bA_{k_{i_\ell}} \bA_{k_{j_{\si(\ell)}}}^*\big)}+\oo1 .
\eeqy
Indeed, similarly from above, we use classical H\"older's inequality to state that the difference between
$$
N^{-n} \Ec{\prod_{\ell=1}^n \tr\big(\bA_{k_{i_\ell}} \bA_{k_{j_{\si(\ell)}}}^*\big)}
$$
and
$$
N^{-(n-1)}\Ec{\prod_{i=1}^{n-1} \tr(\bA_{k_{i_\ell}} \bA_{k_{j_{\si(\ell)}}}^*\big)}  \Ec{\ff N \tr(\bA_{k_{i_n}} \bA_{k_{j_{\si(n)}}}^*)}
$$
is lower than
$$
 \prod_{i=1}^{n-1} \Ec{ \big(\ff N\tr\big(\bA_{k_{i_\ell}} \bA_{k_{j_{\si(\ell)}}}^*\big)\big)^{2(n-1)}}^{\ff{2(n-1)}} \var\big(\ff N \tr(\bA_{k_{i_n}} \bA_{k_{j_{\si(n)}}}^*)\big), 
$$
which tends to $0$ thanks to the non-commutative H\"older's inequality and the fact that $\ff N \tr(\bA_n \bA_{\si(n)}^*)$ converges in probability to a constant. We conclude the proof \eqref{115420112015} with a simple induction. Once we have \eqref{115420112015}, we can conclude using the Wick Formula.
\hfill$\square$

\subsection{Proofs of Theorem \ref{theo08521412020154} and Theorem \ref{thbelow0911141210002015}.} \label{section141212015fga}

In this section, we will directly apply \cite[Theorem 2.3 and Theorem 2.10]{Rochet} in order to prove both Theorems \ref{theo08521412020154}  and \ref{thbelow0911141210002015}. To do so, we only need to prove that the {Gaussian elliptic ensemble} satisfies the assumptions of \cite[Theorem 2.3 and Theorem 2.10]{Rochet}. This is the purpose of the following proposition.
\begin{propo}\la{propo17121518h40}
Let $\bX_N := \ff{\sqrt N} \bY_N$ where $\bY_N$ is an $N \ti N$  {Gaussian elliptic matrix} of parameter $\rho$. Then, as $N\to\infty$: 
\bgt
\item[(i)] $ \ds  \norm{\bX_N}$ converges in probability to $1 + |\rho|$; 
\item[(ii)] for any $\delta>0$, as $N$ goes to infinity, we have the convergence  in probability
$$
\sup_{|z|> 1 + |\rho| + \delta} \max_{1 \leq i,j \leq 2r} \lf| \be_i^* \inve{z\bI - \bX_N} \be_j - \delta_{ij}m(z) \ri| \ \longrightarrow \ 0,
$$ 
where $\ds m(z) \ := \ \int \ff{z-w} \mu_\rho(\mathrm{d}w)$;
\item[(iii)] for any $z$ such that $|z|>1+|\rho|+\ep$,  we have the convergence  in probability

$$
\sqrt{N}\Big( \ff N \tr \inve{z - \bX_N} - m(z)\Big) \ \tto \ 0;
$$
\item[(iv)]  the finite marginals of random process $$\ds\lf( \sqrt{N}\big( \be_i^*\inve{z - \bX_N}\be_j - \delta_{ij}\ff N \tr \inve{z - \bX_N}\big)\ri)_{\substack{ |z|> 1+ |\rho| + \ep\\ 1 \leq i,j \leq 2r}}$$ converge to the ones of the complex centered Gaussian process $$\big( \mathcal{G}_{i,j,z}\big)_{\substack{ |z|> 1+ |\rho| + \ep\\ 1 \leq i,j \leq 2r}}$$ satisfying 
\beq
\E\left[\mathcal{G}_{i,j,z} {\mathcal{G}_{{i'},{j'},{z'}}}\right] \; = \; \delta_{i j'} \delta_{i' j} \lf(\int \ff{(z -w)(z'-w)}  \mu_\rho(\mathrm{d}w) - \int \ff{z -w}  \mu_\rho(\mathrm{d}w)\int \ff{z'-w}  \mu_\rho(\mathrm{d}w)\ri) ,  \\
\E\left[\mathcal{G}_{i,j,z} \ol{\mathcal{G}_{{i'},{j'},{z'}}}\right] \; = \; \delta_{i i'} \delta_{j j' } \lf(\int \ff{(z -w)(\ol{z'}-\ol{w})}  \mu_\rho(\mathrm{d}w) - \int \ff{z -w}  \mu_\rho(\mathrm{d}w)\ol{\int \ff{z'-w}  \mu_\rho(\mathrm{d}w)}\ri) ;\eeq   \\
\item[(v)] for any $p \geq 1$, any $1\le i,j \leq 2r$ and any $|z|>1+|\rho|+\ep$, the sequence
$$
\sqrt{N}\big( \be_i^*\big(z - \bX_N\big)^{-p}\be_j - \delta_{ij}\ff N \tr \big(z - \bX_N\big)^{-p}\big)
$$
is tight. 
\ent
\end{propo}

\begin{rem}
One should be careful about the fact that our $m(z)$ is not the same from \cite[Lemma 4.3]{or13} but the opposite. Moreover, for any $|\tta| > 1$, we still have (see \cite[Eq. (5.2) and (5.3)]{or13})
$$
m(z) \ = \ \ff \tta \ \ \Longleftrightarrow \ \ z  \ = \ \tta + \f{\rho}{\tta}, 
$$
so that it is easy to compute $\E\left[\mathcal{G}_{i,j,z} {\mathcal{G}_{{i'},{j'},{z'}}}\right]$ in (iv) for $z = \tta + \f{\rho}{\tta} $ and $z' = \tta' + \f{\rho}{\tta'} $, indeed
\beq
&&\int \ff{(z-w)(z'-w)}  \mu_\rho(\mathrm{d}w) - \int \ff{z-w}  \mu_\rho(\mathrm{d}w) \int \ff{z'-w}  \mu_\rho(\mathrm{d}w)  \\ & = & -\f{m(z)-m(z')}{z - {z'}} - m(z)m({z'}) \ = \ \ff{\tta\tta' - \rho} - \ff{\tta \tta'}.
\eeq
Also, for any $|z|> 2 \sqrt{|\rho|}$, it might be useful to write 
$$
m(z) \ = \ \sum_{k \geq 0} \rho^k \Cat(k) z^{-2k-1},
$$
where $\Cat(k) = \ff{k+1}\binom{2k}{k}$ is the $k$-th Catalan number. 
\end{rem}

\noindent{\bf Proof of Proposition \ref{propo17121518h40}.}
First, (i) is an adaptation of \cite[Th. 2.2]{or13} to the complexe case, whose proof goes exactly along the same lines (as we work with Gaussian entries, the proof is even easier). It implies that, with a probability tending one, for any  fixed $|z| > 1 + |\rho|$, one can write
$$
\inve{z - \bX_N} \ = \ \sum_{k \geq 0} z^{-k-1} \bX^k_N.
$$
 Moreover, if we apply the Theorem \ref{Th_Main} with $$(\bM_\ell) = \big(\sqrt N \bE_{ji}\big)_{1 \leq i,j \leq 2r}\quad\trm{ and }\quad (\bA_k) = \big( \inve{z-\bX_N}-\ff{N}\Tr \inve{z-\bX_N}\bI\big)_{|z|>1+|\rho|+\ep}$$
(since $\bX_N$ is invariant in distribution by unitary conjugation, so is  $\lf( z - \bX_N\ri)^{-p}$ for any $p \geq 1$), we easily obtain (iv). The same for    (v) by changing the exponent $-1$ into $-p$.  At last, we just need to prove (ii) and (iii).


\noindent{\bf Proof of (iii).}
First of all, let us  write, for any $\eta >0$,
\beq
&&\pro\lf( \sqrt N \big|\ff N \tr \inve{z - \bX_N} - m(z)\big| > {\eta}\ri) \\ & \leq & \f{4N}{\eta^2 } \lf(\E \lf| \ff N \tr \inve{z - \bX_N} - \E \ff N \tr \inve{z - \bX_N}\ri|^2 + \lf| \E \ff N \tr \inve{z - \bX_N} - m(z)\ri|^2 \ri),
\eeq
which means that we just have to prove that 
\beqy \la{08321712120157l72015}
\lf| \E \ff N \tr \inve{z - \bX_N} - m(z)\ri|^2 & = & \oo{\ff N},
\eeqy
and 
\beqy \label{2412235178422015}
\E \lf| \ff N \tr \inve{z - \bX_N} - \E  \ff N \tr \inve{z - \bX_N}\ri|^2  & =  & \oo{\ff N}.
\eeqy

\noindent\underline{Proof of \eqref{08321712120157l72015}:} We know from \cite[Theorem 1.1]{BAI} that \eqref{08321712120157l72015} woud be true had $\bX_N$ been a Gaussian Wigner matrix instead of an elliptic one. Here, the idea of the proof is to use the fact that the Stieltjes transform of the semicircular law of variance $\si^2 = \rho$ is equal to $m(z)$ outside the ellipse $\mathcal{E}_\rho$ when $\rho >0$.
First, we shall suppose that $\rho \geq 0$ up to changing $\bX_N$ into $\ii\bX_N$. 
For any $i \neq j$, we have 
\beqy \label{1426201150914}
\E x_{ii}^2 \ = \  \E x_{ij}{x_{ji}} \ = \ \f{\rho}{N}, & \text{and} & \E x_{ij}^2 \ = \ 0.
\eeqy
 One can notice that if  $\bW_N$ is  a real symmetric Gaussian matrix of variance $\rho$ with iid entries such that, for any $i \neq j$,
\beqy \label{1426201150914bis}
 \E w_{ii}^2 \  = \ \E w_{ij}{w_{ji}} \ = \ \E w_{ij}^2 \ = \ \f{\rho}{N},
\eeqy
 then
we have, by the Wick formula applied to the expansion of the traces, 
\bgt
\ite[$\bullet$] $\Ec{\tr  \bW_N^{k}}\ge 0$ and $\Ec{\tr \bX_N^{k}}\ge 0$,
\ite[$\bullet$] $\Ec{\tr  \bW_N^{k}} \geq \Ec{\tr \bX_N^{k}}$ since there are more non-zero terms for $\bW_N$ than for $\bX_N$.
\ent
Also, we know that, for any $z$ such that $|z| > 1 + \rho+\ep$, $$\E\ff N \tr\inve{z-\bW_N} \ = \  \sum_{k \geq 0} z^{-k-1} \Ec{\ff N \tr \bW_N^{k}} \ \text{   and   } \ \E\ff N \tr \inve{z-\bX_N} \ = \ \sum_{k \geq 0} z^{-k-1} \Ec{\ff N \tr \bX_N^{k}}$$ converge to the same limit $$m(z) \ = \ \sum_{k \geq 0} \Cat(k) \rho^k z^{-2k-1}, \ \ \ \text{ where } \ \Cat(k)\ \text{ is the } k \text{-th Catalan number.} $$  Moreover, 
by the Wick formula again, for $\mathcal{P}_2(2k)$ (resp. ${NC}_2(2k)$) the set of pairings (resp. non crossing pairings) of $\{1, \ld,2k\}$,  
\beq
\Ec{\tr \bX_N^{2k}} & = & \sum_{1 \leq i_1,\ld,i_{2k} \leq N} \sum_{\pi \in \mathcal{P}_2(2k)} \prod_{\{s,t\} \in \pi} \Ec{x_{i_s i_{s+1}} x_{i_t i_{t+1}}} \\
& = & \sum_{\pi \in \mathcal{P}_2(2k)} \sum_{1 \leq i_1,\ld,i_{2k} \leq N}  \prod_{\{s,t\} \in \pi} \Ec{x_{i_s i_{s+1}} x_{i_t i_{t+1}}} \eeq 
Note that using the Dyck path interpretation of ${NC}_2(2k)$ (see e.g. \cite{ns06}), one can easily see that in the previous sum, the term associated to each  $\pi \in NC_2(2k)$ is precisely $\rho^k$. Hence as the cardinality of $NC_2(2k)$ is $\Cat(k)$ (see  \cite{ns06} again) and each $ \Ec{x_{i_s i_{s+1}} x_{i_t i_{t+1}}} $ is non negative, we have 
\beq
\Ec{\ff N \tr \bX_N^{2k}}& \ge & \Cat(k) \rho^k  .
\eeq
At last, we know from e.g.  \cite[Theorem 1.1]{BAI} that, for any $z$ such that $ |z|>1+\rho+\ep$, we have 
\beq
\E\ff N \tr \inve{z-\bW_N} - m(z) & = & \oo{\ff{\sqrt N}},
\eeq 
so that, to conclude, it suffices to write
\beq
\lf|\E\ff N \tr \inve{z-\bX_N} - m(z)\ri| & \leq & \sum_{k \geq 0} \big( \E \ff N \tr \bX_N^{2k} - \Cat(k)\rho^k\big)|z|^{-2k-1} \\
& \leq & \sum_{k \geq 0} \big( \E \ff N \tr \bW_N^{2k} - \Cat(k)\rho^k\big)|z|^{-2k-1} \\  
& = & \E\ff N \tr \inve{|z|-\bW_N} - m(|z|) \ = \ \oo{\ff{\sqrt N}}. \\ \eeq
\underline{Proof of \eqref{2412235178422015}:} We apply the same idea but this time $\bW_N$ is  a real symmetric Gaussian matrix of variance $\rho$ with iid entries, such that, for any $i \neq j$,
\beqy \label{1426201150914ter}
 \E w_{ii}^2 \  = \ \ff N &;& \E w_{ij}{w_{ji}} \ = \ \E w_{ij}^2 \ = \ \f{\rho}{N}.
\eeqy 
From e.g. \cite[Theorem 1.1]{BAI}, we know that, for all $|z| > 1 + |\rho| + \ep$,
$$
\E \lf| \ff N \tr \inve{z - \bW_N} - \E  \ff N \tr \inve{z - \bW_N}\ri|^2  \ =  \ \oo{\ff N}. 
$$
Moreover, we can write
\beq
&&\E \lf| \ff N \tr \inve{z - \bX_N} - \E  \ff N \tr \inve{z - \bX_N}\ri|^2 \\ & =  & \E \lf| \ff N \tr \inve{z - \bX_N}\ri|^2 - \lf| \E \ff N \tr \inve{z - \bX_N} \ri|^2 \\ 
& = & \ff{ N^2 }\sum_{k,\ell \geq 0} z^{-k-1} \ol{z}^{-\ell-1} \lf(\Ec{\tr \bX_N^k \ol{\tr \bX_N^\ell}  } -  \E\tr \bX_N^k \E\ol{\tr \bX_N^\ell}  \ri).
\eeq
By  the Wick formula, we   see  that, for all $k,\ell$,
$$
0 \ \leq \ \Ec{\tr \bX_N^k \ol{\tr \bX_N^\ell}  } -  \E\tr \bX_N^k \E\ol{\tr \bX_N^\ell}  \ \leq \ \Ec{\tr \bW_N^k \ol{\tr \bW_N^\ell}  } -  \E\tr \bW_N^k \E\ol{\tr \bW_N^\ell}  .
$$
Indeed, 
\beqy\nonumber
\Ec{\tr \bX_N^k \ol{\tr \bX_N^\ell}  } & = & \sum_{\substack{1 \leq i_1,\ld, i_k \leq N \\ 1 \leq i'_{k+1},\ld,i'_{k+\ell} \leq N}} \Ec{x_{i_1 i_2} \cdots x_{i_k i_1} \ol{x_{i'_{k+1} i'_{k+2}}} \cdots \ol{x_{i'_{k+\ell} i'_{k+1}}}}\\ \label{25121519h05}
& = & \sum_{\substack{1 \leq i_1,\ld, i_k \leq N \\ 1 \leq i'_{k+1},\ld,i'_{k+\ell} \leq N}} \sum_{\pi \in \mathcal{P}_2(k+\ell)} \prod_{\{a,b\} \in \pi} \Ec{x_{a} x_{b}} 
\eeqy
where $x_a = \begin{cases} x_{i_a i_{a+1}} & \text{ if } 1 \leq a \leq k-1 \\
                           x_{i_k i_{1}} &  \text{ if } a= k \\ \ol{x_{i'_a i'_{a+1}}} & \text{ if } k+1 \leq a \leq k+\ell-1 \\
                           \ol{x_{i'_{k+\ell} i'_{k+1}}} &  \text{ if } a= k +\ell. \\            \end{cases}$ \\
               We have also
\beq
\E\tr \bX_N^k \E\ol{\tr \bX_N^\ell} & = & \sum_{\substack{1 \leq i_1,\ld, i_k \leq N \\ 1 \leq i'_{k+1},\ld,i'_{k+\ell} \leq N}} \sum_{\substack{\pi \in \mathcal{P}_2(k) \\ \mu \in \mathcal{P}_2(\ell)}} \prod_{\{a,b\} \in \pi} \Ec{x_{a} x_{b}} \prod_{\{c,d\} \in \mu} \Ec{x_{c} x_{d}},\\
\eeq
which is a subsum of \eqref{25121519h05}. Hence,   as all the $\Ec{x_a x_b}$'s are non negative (see \eqref{1426201150914}), we conclude that 
$$
\Ec{\tr \bX_N^k \ol{\tr \bX_N^\ell}  } \ \geq \ \E\tr \bX_N^k \E\ol{\tr \bX_N^\ell}.
$$ 
Since, for all $a$ and $b$, $\Ec{x_a x_b} \leq \Ec{w_a w_b}$ (see \eqref{1426201150914ter}), we deduce that
$$
\Ec{\tr \bX_N^k \ol{\tr \bX_N^\ell}  } -  \E\tr \bX_N^k \E\ol{\tr \bX_N^\ell}  \ \leq \ \Ec{\tr \bW_N^k \ol{\tr \bW_N^\ell}  } -  \E\tr \bW_N^k \E\ol{\tr \bW_N^\ell}
$$
At last, we can write 
\beq
&&\E \lf| \ff N \tr \inve{z - \bX_N} - \E  \ff N \tr \inve{z - \bX_N}\ri|^2 \\ 
& \leq  & \ff{ N^2 }\sum_{k,\ell \geq 0} |z|^{-k-1} |z|^{-\ell-1} \lf(\Ec{\tr \bX_N^k \ol{\tr \bX_N^\ell}  } -  \E\tr \bX_N^k \E\ol{\tr \bX_N^\ell}  \ri)\\
& \leq & \ff{ N^2 }\sum_{k,\ell \geq 0} |z|^{-k-1} |z|^{-\ell-1} \lf(\Ec{\tr \bW_N^k \ol{\tr \bW_N^\ell}  } -  \E\tr \bW_N^k \E\ol{\tr \bW_N^\ell}  \ri)\\
& = & \E \lf| \ff N \tr \inve{|z| - \bW_N} - \E  \ff N \tr \inve{|z| - \bW_N}\ri|^2 \ = \ \oo{\ff N}.
\eeq
\hfill $\square$

\noindent{\bf Proof of (ii).}
Let $\eta>0$ and let $i,j$ be two integers lower than $2r$. Since $\norm{\bX_N}$ is bounded, we know that $\norm{\inve{z-\bX_N}}$ goes to $0$ when $|z|\to \infty$, as the function $m(z)$, so that, we know there is a positive constant $M$ such that
$$
\pro\big( \hspace{-2mm}\sup_{|z|> 1 + \rho+\ep} \hspace{-4mm} \big| \be_i^*\inve{z - \bX_N}\be_j - \delta_{ij}m(z) \big| > \eta \big)  =  \pro\big(\hspace{-2mm} \sup_{  1 + \rho+\ep <|z|<M} \hspace{-5mm} \big| \be_i^*\inve{z - \bX_N}\be_j - \delta_{ij}m(z) \big| > \eta \big) + \oo1.
$$
Then, for any $\eta'>0$,  the compact set $K = \{1 + \rho + \ep\leq |z| \leq M \}$ admits a $\eta'$-net that we denote by $S_{\eta'}$, which is a finite set of $K$ such that
$$
\forall z \in K, \ \exists z' \in S_{\eta'}, \ \ \ |z-z'| < \eta',
$$
so that, using the uniform boundedness of the derivative  of $m(z)$ and $\be_i^*\inve{z - \bX_N}\be_j$ on $K$, we have for a small enough $\eta'$
\beq
\pro\big( \sup_{ z \in K}  \big| \be_i^*\inve{z - \bX}\be_j - \delta_{ij}m(z) \big| > \eta \big) & = & \pro\big( \bigcup_{z \in S_{\eta'}} \lf\{ \big| \be_i^*\inve{z - \bX}\be_j - \delta_{ij}m(z) \big| > \eta/2 \ri\}\big)
\eeq
At last, we write, for any $z \in S_{\eta'}$,
\beq
&&\pro\big(  \lf| \be_i^*\inve{z - \bX}\be_j - \delta_{ij}m(z) \ri| > \eta/2 \big)\\& \leq & \pro\big(  \lf| \be_i^*\inve{z - \bX}\be_j - \delta_{ij}\ff N \tr \inve{z - \bX} \ri| > \eta/4 \big) +\pro\big(  \delta_{ij}\lf| \ff N \tr \inve{z - \bX}-m(z) \ri| > \eta/4 \big).
\eeq
The first term vanishes thanks to Theorem \ref{Th_Main} with $\bM = \sqrt N \bE_{ji}$ and the second one vanishes by (iii).
\hfill $\square$

\section*{Appendix: a matrix inequality}

\begin{lemme}\label{lem101206102015}
For any $k\ge 2$ and any Hermitian matrix $\bH$, 
$$
\lf|\tr \bH^{k}\ri| \ \leq \ \big(\tr \bH^2\big)^{k/2} .
$$
More generally, for any family of Hermitian matrices $\bH_1,\ldots, \bH_k$, 
$$
\lf| \tr \big( \bH_1 \cdots \bH_k\big) \ri| \ \leq \ \prod_{i=1}^k \sqrt{\tr \bH_i^2}
$$
\end{lemme}
\bpr
We know that, for any non negative Hermitian matrices $\bA$ and $\bB$, one has
$$
\tr \bA \bB \leq \tr \bA \tr \bB
$$ 
so that, for any $p\ge 1$, 
$$
\tr \bH^{2p} \leq \lf(\tr \bH^2\ri)^p
$$
also, 
$$
\tr \bH^{2p+1} \leq \sqrt{\tr \bH^2} \sqrt{\tr\bH^{4p}} \leq \sqrt{\tr \bH^2} \sqrt{\big(\tr\bH^{2}\big)^{2p}} \ = \ \big(\tr \bH^2\big)^{(2p+1)/2}. 
$$
Then, using the non-commutative H\"older's inequality (see \cite[A.3]{agz}), we deduce that
$$
\lf|\tr \big( \bH_1 \cdots \bH_k\big) \ri|\leq \prod_{i=1}^k \big( \tr |\bH_i|^k \big)^{1/k} \leq \prod_{i=1}^k \Big(\big( \tr \bH_i^2 \big)^{k/2}\Big)^{1/k}.
$$
\epr

\end{document}